\documentclass[12pt,amssymb,verbatim]{amsart}

\usepackage{amsfonts,latexsym,geometry,color,hyperref,ulem}
\usepackage{amssymb,graphics,graphicx}
\usepackage{amsfonts,latexsym,geometry}
\usepackage{amssymb,graphics,graphicx}
\usepackage{bbm}
\usepackage{amscd}
\usepackage{amssymb}
\usepackage{amsthm}
\usepackage[utf8]{inputenc} 
\usepackage[T1]{fontenc}
\usepackage{amsmath}
\usepackage{amsfonts}
\usepackage{amssymb}
\usepackage{graphicx} 
\usepackage{enumerate}
\newcommand\supp{\mathrm{supp}}
\newtheorem{theoreme}{Theorem}[section] %
\newtheorem{proposition}[theoreme]{Proposition} %
\newtheorem{corollary}[theoreme]{Corollary} %
\newtheorem{lemme}[theoreme]{Lemma} %
\newtheorem{definition}{Definition}[section] %

\newcommand\sk{\smallskip}

\newcommand\R{\mathbb{R}} 
 
\newcommand\dimm{\underline{\dim}_H}

\RequirePackage{yhmath} 
\renewcommand\widering[1]{\ring{#1}}

\author{E.Daviaud, ULiège}

\begin{document}
\title[Slow polunomial mixing]{Slow polynomial mixing, dynamical Borel Cantelli-Lemma and Hausdorff dimension of dynamical Diophantine sets}
\maketitle
\begin{abstract}
In \cite{Edergo}, it is established that if an ergodic system  $(T,\mu)$ of $\mathbb{R}^d$,  where $\mu$ is an exact-dimensional measure, satisfies a certain mixing property, then  the dynamical Borel-Cantelli lemma holds at every point and for every $\delta\geq \frac{1}{\dim \mu},$ for $\mu$-almost every $x_0$, $\dim_H \left\{x:\vert\vert T^n(x_0)-x\vert\vert_{\infty}\leq \frac{1}{n^{\delta}}\text{ i.o. }\right\}=\frac{1}{\delta}.$  
In the present article, improving on an example of Galatolo-Rousseau-Saussol, we show that the mixing condition provided in \cite{Edergo} is optimal in a satisfying sense for both the dynamical Borel-Cantelli lemma and to get the above estimate regarding the dimension of the classical dynamical Diophantine sets.
\end{abstract}
\section{Introduction}
The Diophantine approximation takes his roots in number theory, with the work of Dirichlet, who  defined and proved the first result regarding the notion of speed of approximation of a real number by rationals. More precisely, Dirichlet proved that for any $x\in \mathbb{R}\setminus\mathbb{Q},$ there exists infinitely many $(p,q)$ with $p\wedge q=1$ such that $$\vert x-\frac{p}{q}\vert\leq \frac{1}{q^{2}}.$$
His result was later on refined by Jarnik and Besicovitch by establishing that for every $\delta \geq 1,$ $$\dim_H \left\{x: \vert x -\frac{p}{q}\vert \leq\frac{1}{q^{2\delta}}\text{ i.o. }p\wedge q=1\right\}=\frac{1}{\delta},$$
where $``$ i.o.$"$ (infinitely often) means that the inequality holds for infinitely many $p,q$. Finally recently, by combining an analog of Borel-Cantelli lemma proved by Koukoulopoulos and Maynard \cite{Maynkou} and the so-called mass transference principle from Beresnevitch and Velani \cite{BV} or Jaffard's ubiquity theorem \cite{Jaff-prescription}, it was proved that, given $\psi:\mathbb{N}\to\mathbb{R}_+,$ $$\dim_H \left\{x: \vert x -\frac{p}{q}\vert \leq\psi(q)\text{ i.o. }p\wedge q=1\right\}=\min\left\{1,s_{\psi}\right\},$$
where, denoting $\phi$ the Euler mapping, $$s_{\psi}=\inf\left\{s: \sum_{q\geq 1}\phi(q)\psi(q)^s<+\infty\right\}.$$
Similar questions have also raised many interests in dynamical systems, as analogs of rational approximation by orbits are very natural.  Given a measurable mapping $T:\mathbb{R}^d \to \mathbb{R}^d$ (where $d\in\mathbb{N}$) a $T$-ergodic (probability) measure $\mu$, a point $x_0 \in\mathbb{R}^d$ and a non-increasing sequence of radii $(r_n)_{n\in\mathbb{N}}\in \mathbb{R}_+ ^{\mathbb{N}}$, two particular problems led to many interesting developments, which consists in studying, from the measure or Hausdorff dimension point of view the sets $$\begin{cases}\left\{x: \ T^{n}(x)\in B(x_0,r_n)\text{ i.o.}\right\} \\ \left\{x : \ \vert\vert T^n (x_0)-x\vert\vert_{\infty}\leq r_n \text{ i.o.}\right\}.\end{cases} $$
 Regarding the first problem, a very natural question consists in investigating under which condition on the ergodic system $(T,\mu)$ one obtains an analog of the Borel-Cantelli lemma at the point $x_0,$ i.e., to determine whether it holds that 
 \begin{equation}
 \label{dichoBorelCant}
 \begin{cases}\mu\Big(\left\{x: \ T^{n}(x)\in B(x_0,r_n)\text{ i.o. }\right\}\Big)=1 \text{ if }\sum_{n\geq 0} \mu\Big(B(x_0,r_n)\Big)=+\infty \\ \mu\Big(\left\{x: \ T^{n}(x)\in B(x_0,r_n)\text{ i.o. }\right\}\Big)=0 \text{ if }\sum_{n\geq 0} \mu\Big(B(x_0,r_n)\Big)<+\infty.\end{cases} \end{equation}
 The answer regarding this question for certain conformal Gibbs measure where established by Chernov and Kleinbock in \cite{ChernovKleinbock} and, later on, Galatolo \cite{Galadim} established that this dichotomy holds at any point provided that $(T,\mu)$ is exponentially mixing with respect to Lipschitz observable and $\mu$ is exact-dimensional (see Definition \ref{dimmu}). In \cite[Proposition 41]{Saussogalato} Galatolo-Saussol-Rousseau provided an example (a class of examples to be more precise) of ergodic systems $(T,\mu)$, where $T:\mathbb{T}^3 \to\mathbb{T}^3$ ($\mathbb{T}$ denoting the one dimensional torus) and $\mu=\mathcal{L}^3$ is the Lebesgue measure which are polynomially mixing with respect to Lipschitz-observable but for which at any point $x_0$, 
 $$\mu\Big(\left\{x: \ T^{n}(x)\in B(x_0,\frac{1}{n^{\frac{1}{3}}})\text{ i.o. }\right\}\Big)=0,$$
 which shows that at every $x_0$, the Borel-Cantelli dichotomy \eqref{dichoBorelCant} does not hold.
 
In \cite{Edergo}, following \cite{LLVZ}, was introduced the $(\gamma,\phi)$-mixing property, where   $\gamma \geq 1$ and $\phi:\mathbb{N}\to \mathbb{R}_+$ is a non increasing mapping. An ergodic system is said to be $(\gamma,\phi)$-mixing if for every balls $A,B$\footnote{It is worth mentioning that, by outer regularity, it is easily proved that the ball $A$ can be replaced by any Borel set in \eqref{defgammaphimix}, but not the ball $B$. Indeed, if $\mu$ is an ergodic measure satisfying this properties for every pairs of Borel sets $A,B$, then it is a small exercise to show that $\mu =\frac{1}{N}\sum_{0\leq k\leq N}\delta_{T^k(y)}$ for some periodic orbit $\left\{T^k(y)\right\}_{0\leq k\leq N}$.}, for every $n\in\mathbb{N}$ one has 
\begin{equation}
\label{defgammaphimix}
\mu\Big(T^{-n} (A)\cap B \Big)\leq \gamma \mu(A)\mu(B)+\phi(n)\mu(A).
\end{equation}
The advantage of \eqref{defgammaphimix} is that it is often relatively easy to verify whether an ergodic $(T,\mu)$ satisfies it or not. Moreover, should $(T,\mu)$ satisfy \eqref{defgammaphimix} with $\phi(n)=Cn^{1+\varepsilon}$, for some $\varepsilon>0,$ combining Chung-Erdös inequality and Birkhoff's pointwise ergodic theorem (see \cite[Lemma 6.2]{Edergo}), one proves that \eqref{dichoBorelCant} holds at any point and for every non-increasing sequence of radii $(r_n)_{n\in\mathbb{N}}$ and it was established in \cite[Theorem 3.2]{Edergo} that, if in addition $\mu$ is exact-dimensional, then for every $\delta\geq \frac{1}{\dim_H \mu},$ for $\mu$-almost every $x_0$, one has
\begin{equation}
\label{DimDioDyna}
\dim_H \left\{x : \ \vert\vert T^n (x_0)-x\vert\vert_{\infty}\leq \frac{1}{n^\delta}\text{ i.o.}\right\}=\frac{1}{\delta}.
\end{equation}
In this article we establish that there exists $0< s\leq 1$ and $C>0$ such that the example $(T,\mathcal{L}^3)$ provided by Galatolo-Saussol-Rousseau in \cite{Saussogalato} is $(9, n\mapsto \frac{C}{n^s})$-mixing and that, for Lebesgue-almost every $x_0$, $$\dim_H \left\{x : \ \vert\vert T^n (x_0)-x\vert\vert_{\infty}\leq \frac{1}{n^{\frac{1}{3}}}\text{ i.o.}\right\}<3.$$
This shows that, given an ergodic system $(T,\mu)$ with  exact-dimensional $\mu$, both \eqref{dichoBorelCant} and \eqref{DimDioDyna} holds as soon as $(T,\mu)$ is $(\gamma,n\mapsto \frac{C}{n^{1+\varepsilon}})$-mixing and needs not hold if $(T,\mu)$ is $(\gamma,n\mapsto \frac{C}{n^{1-\varepsilon}})$-mixing.
\section{Statement of the main result}

Let us start with some notations 

 Let $d$ $\in\mathbb{N}$. For $x\in\mathbb{R}^{d}$, $r>0$,  $B(x,r)$ stands for the closed ball of ($\mathbb{R}^{d}$,$\parallel \ \ \parallel_{\infty}$) of center $x$ and radius $r$. 
 Given a ball $B$, $\vert B\vert$ stands for the diameter of $B$. For $t\geq 0$, $\delta\in\mathbb{R}$ and $B=B(x,r)$,   $t B$ stands for $B(x,t r)$, i.e. the ball with same center as $B$ and radius multiplied by $t$,   and the  $\delta$-contracted  ball $B^{\delta}$ is  defined by $B^{\delta}=B(x ,r^{\delta})$.
\smallskip

Given a set $E\subset \mathbb{R}^d$, $\widering{E}$ stands for the  interior of the set $E$, $\overline{E}$ its  closure and $\partial E =\overline{E}\setminus \widering{E}$ its boundary. If $E$ is a Borel subset of $\R^d$, its Borel $\sigma$-algebra is denoted by $\mathcal B(E)$.
\smallskip

Given a topological space $X$, the Borel $\sigma$-algebra of $X$ is denoted $\mathcal{B}(X)$ and the space of probability measure on $\mathcal{B}(X)$ is denoted $\mathcal{M}(X).$ 

\sk
 The $d$-dimensional Lebesgue measure on $(\mathbb R^d,\mathcal{B}(\mathbb{R}^d))$ is denoted by 
$\mathcal{L}^d$.
\smallskip

For $\mu \in\mathcal{M}(\R^d)$,   $\supp(\mu)=\left\{x\in \mathbb{R}^d: \ \forall r>0, \ \mu(B(x,r))>0\right\}$ is the topological support of $\mu$.
\smallskip

 Given $E\subset \mathbb{R}^d$, $\dim_{H}(E)$ and $\dim_{P}(E)$ denote respectively  the Hausdorff   and the packing dimension of $E$.
\smallskip

Given a set $S$, $\chi_S$ denotes the indicator function of $S$, i.e., $\chi_S(x)=1$ if $x\in S$ and $\chi_S(x)=0$ otherwise.

 \smallskip

 Given $n\in\mathbb{N},$ $\mathcal{D}_n$ denotes the set of dyadic cubes of generation $n$ and $\mathcal{D}$ the set of all dyadic cubes, i.e.
 \begin{equation*}
\mathcal{D}_{n}=\left\{2^{-n}(k_1,...,k_d)+2^{-n}[0,1)^d, \ (k_1,...,k_d)\in\mathbb{Z}^d\right\}\text{ and }\mathcal{D}=\bigcup_{n\geq 0}\mathcal{D}_{n}.
\end{equation*}

\smallskip

An ergodic system will always refer to a couple $(T,\mu)$ where $\mu$ is measurable and $\mu$ is $T$-ergodic.

\subsection{Recall on geometric measure theory}

\begin{definition}
\label{hausgau}
Let $\zeta :\mathbb{R}^{+}\mapsto\mathbb{R}^+$. Suppose that $\zeta$ is increasing in a neighborhood of $0$ and $\zeta (0)=0$. The  Hausdorff outer measure at scale $t\in(0,+\infty]$ associated with the gauge $\zeta$ of a set $E$ is defined by 
\begin{equation}
\label{gaug}
\mathcal{H}^{\zeta}_t (E)=\inf \left\{\sum_{n\in\mathbb{N}}\zeta (\vert B_n\vert) : \, \vert B_n \vert \leq t, \ B_n \text{ closed ball and } E\subset \bigcup_{n\in \mathbb{N}}B_n\right\}.
\end{equation}
The Hausdorff measure associated with $\zeta$ of a set $E$ is defined by 
\begin{equation}
\mathcal{H}^{\zeta} (E)=\lim_{t\to 0^+}\mathcal{H}^{\zeta}_t (E).
\end{equation}
\end{definition}

For $t\in (0,+\infty]$, $s\geq 0$ and $\zeta:x\mapsto x^s$, one simply uses the usual notation $\mathcal{H}^{\zeta}_t (E)=\mathcal{H}^{s}_t (E)$ and $\mathcal{H}^{\zeta} (E)=\mathcal{H}^{s} (E)$, and these measures are called $s$-dimensional Hausdorff outer measure at scale $t\in(0,+\infty]$ and  $s$-dimensional Hausdorff measure respectively. Thus, 
\begin{equation}
\label{hcont}
\mathcal{H}^{s}_{t}(E)=\inf \left\{\sum_{n\in\mathbb{N}}\vert B_n\vert^s : \, \vert B_n \vert \leq t, \ B_n \text{  closed ball and } E\subset \bigcup_{n\in \mathbb{N}}B_n\right\}. 
\end{equation}

\begin{definition} 
\label{dim}
Let $\mu\in\mathcal{M}(\mathbb{R}^d)$.  
For $x\in \supp(\mu)$, the lower and upper  local dimensions of $\mu$ at $x$ are  defined as
\begin{align*}
\underline\dim_{{\rm loc}}(\mu,x)=\liminf_{r\rightarrow 0^{+}}\frac{\log(\mu(B(x,r)))}{\log(r)}
 \mbox{ and } \ \    \overline\dim_{{\rm loc}}(\mu,x)=\limsup_{r\rightarrow 0^{+}}\frac{\log (\mu(B(x,r)))}{\log(r)}.
 \end{align*}
Then, the lower and upper Hausdorff dimensions of $\mu$  are defined by 
\begin{equation}
\label{dimmu}
\dimm(\mu)={\mathrm{ess\,inf}}_{\mu}(\underline\dim_{{\rm loc}}(\mu,x))  \ \ \mbox{ and } \ \ \overline{\dim}_P (\mu)={\mathrm{ess\,sup}}_{\mu}(\overline\dim_{{\rm loc}}(\mu,x))
\end{equation}
respectively.
\end{definition}

It is known (for more details see \cite{F}) that
\begin{equation*}
\begin{split}
\dimm(\mu)&=\inf\{\dim_{H}(E):\, E\in\mathcal{B}(\mathbb{R}^d),\, \mu(E)>0\} \\
\overline{\dim}_P (\mu)&=\inf\{\dim_P(E):\, E\in\mathcal{B}(\mathbb{R}^d),\, \mu(E)=1\}.
\end{split}
\end{equation*}
When $\underline \dim_H(\mu)=\overline \dim_P(\mu)$, this common value is simply denoted by $\dim(\mu)$ and~$\mu$ is said to be \textit{exact-dimensional}. It is worth mentioning that many measures which are dynamically defined are exact dimensional. 

\subsection{Continued fraction, Diophantine linear type}
 Given $\theta\in [0,1]$ the continued fraction expansion of $\theta$ is defined as the unique sequence of integers $[\theta_n]_{n\geq 1}$ such that $$\theta=\frac{1}{\theta_1+\frac{1}{\theta _2+...}}.$$

Given $n\in\mathbb{N}$, writing $$\frac{1}{\theta_1+\frac{1}{\theta_2+\frac{1}{...+\frac{1}{\theta_n}}}}=\frac{p_n}{q_n}$$
where $p_n \wedge q_n=1$, $\frac{p_n}{q_n}$ is called the $n$th partial quotient of $\theta.$ Moreover it is known that for every $n\in\mathbb{N},$ $$\vert \theta -\frac{p_n}{q_n}\vert\leq \frac{1}{q_n q_{n+1}}.$$
In the rest of the article, the following Diophantine exponent will play an important role.

\begin{definition}
Let $\alpha\in\mathbb{R}^2.$ The Diophantine linear type of $\alpha$ is denoted $\gamma_{\ell}(\alpha)$ and defined as $$\gamma_{\ell}(\alpha)=\inf\left\{\gamma\geq 0 \text{ s.t. } \exists C>0 : \  d_{\mathbb{Z}}(\alpha.k)\geq \frac{C}{\vert\vert k \vert\vert_{\infty}^{\gamma}}\text{ for every  }k\in\mathbb{Z}^2\right\},$$
where $d_{\mathbb{Z}}$ denotes the distance to the nearest integer and $\cdot$ the Euclidien inner-product.
\end{definition}

\subsection{Statement of the main result}

 Let $\alpha\in\mathbb{R}^2$ be a vector with rational independent coordinates and define $S:\mathbb{T}^1\times  \mathbb{T}^2 \to \mathbb{T}^3$ by setting 
 $$ S_{\alpha}(x,t)=(2x,t+\alpha \chi_{[0,\frac{1}{2}[}(x)).$$
It is easily proved that $(S(x,t),\mathcal{L}^3)$ is mixing.

These skew-products where originally introduced in \cite{Saussogalato}, in order to provide a class of dynamical system which are polynomially mixing with respect to Lipschitz observable (provided that $\alpha$ is well chosen) but for which the dynamical Borel-Cantelli fails at every point. 
More precisely, the authors proved first the following:
\begin{proposition}[\cite{Saussogalato}]
\label{alpha12}
There exists $\alpha_1,\alpha_2 \in\mathbb{R}$ such that writing $(\frac{p_n}{q_n})_{n\in\mathbb{N}}$ the partial quotients of $\alpha_1$ and $(\frac{p'_n}{q'_n})_{n\in\mathbb{N}}$ the one of $\alpha_2,$ one has $$\begin{cases}q_n^4 \leq q'_n \leq 4q_n^4 \\ {q'}_{n-1}^4\leq q_n \leq 4{q'}_{n-1}^4 \\
\gamma_{\ell}(\alpha_1,\alpha_2)=16.\end{cases} $$
\end{proposition}

We fix $\alpha=(\alpha_1,\alpha_2),$ where $\alpha_1,\alpha_2$ are given by Proposition \ref{alpha12}. It is proved in \cite{Saussogalato} that $(S_{\alpha},\mathcal{L}^3)$ is polynomially mixing with respect to Lipschitz observable. 

 The authors obtained the following result.
\begin{theoreme}[\cite{Saussogalato}]
\label{contrexBorelCantelli}
Let $\alpha \in\mathbb{R}^2$ be as in Proposition \ref{alpha12}. Then for every $y \in\mathbb{T}^3,$ for every $\delta>\frac{1}{4},$ one has $$\mathcal{L}^{3}\Big(\left\{ x: \ T^{n}(x)\in B\Big(y,\frac{1}{n^{\delta}}\Big)\text{ i.o. }\right\}\Big)=0.$$
\end{theoreme}
Notice that, taking $\delta=\frac{1}{3}$,  $$\sum_{n\geq 1}\mathcal{L}^3\Big(B(y,\frac{1}{n^{\frac{1}{3}}})\Big)=\sum_{n\geq 1}\frac{1}{n}=+\infty$$
whereas $$\mathcal{L}^{3}\Big(\left\{ x: \ T^{n}(x)\in B\Big(y,\frac{1}{n^{\frac{1}{3}}}\Big)\text{ i.o. }\right\}\Big)=0.$$
Thus, this shows the dichotomy provided by \eqref{dichoBorelCant} needs not hold in general. The next result is proved in \cite{Edergo} (see proof of \cite[Lemma 6.2 and Theorem 3.2]{Edergo}).

\begin{theoreme}[\cite{Edergo}]
\label{Thmergo}
Let $(T,\mu)$ be an ergodic system such that $\mu$ is exact-dimensional and $(T,\mu)$ is $(\gamma,\frac{C}{n^{1+\varepsilon}})$-mixing, for some $C,\varepsilon>0$ and $\gamma\geq 1.$ Then:
\begin{itemize}
\item[•] The dichotomy \eqref{dichoBorelCant} holds at every point and for any non-increasing sequence of radii $(r_n)_{n\in\mathbb{N}},$ \medskip
\item[•] for $\mu$-almost every $x_0,$ for every $\delta \geq \frac{1}{\dim_H \mu},$ one has $$\dim_H\left\{x: \ \vert\vert x-T^n(x_0)\vert\vert_{\infty}\text{ i.o.}\right\}=\frac{1}{\delta}.$$
\end{itemize}
\end{theoreme}

Our main result is the following:

\begin{theoreme}
\label{Mainthm}
Let $\alpha\in\mathbb{T}^2$ be as in Proposition \ref{alpha12}. Then:
\begin{itemize}
\item[(1)]There exists $C>0$ and $0 \leq s\leq 1$ such that the ergodic system $(S_{\alpha},\mathcal{L}^3)$ is $(9,n\mapsto \frac{C}{n^{s}})$-mixing (see \eqref{defgammaphimix}),\medskip
\item[(2)] for  every $x_0$, $$\dim_H \left\{x: \ \vert\vert T^n(x_0)-x\vert\vert_{\infty}\leq \frac{1}{n^{\frac{1}{3}}}\text{ i.o.}\right\}<3.$$
\end{itemize}

\end{theoreme}
Theorem \ref{Mainthm} together with Theorem \ref{contrexBorelCantelli} yield the following corollary, confirming that Theorem \ref{Thmergo} is optimal in a strong sens.

\begin{corollary}
There exists an ergodic system $(T,\mu)$ on $\mathbb{T}^3$, with $\mu$ exact-dimensional and such that $(T,\mu)$ is $(\gamma,n\mapsto\frac{C}{n^{1-\varepsilon}})$-mixing for some $C>0,$ $0<\varepsilon<1$ and $\gamma\geq 1$ and for which:
$$\begin{cases} \text{ the dichotomy }\eqref{dichoBorelCant} \text{ fails at every point,} \\ \dim_H \left\{x: \ \vert\vert T^n(x_0)-x\vert\vert_{\infty}\leq \frac{1}{n^{\frac{1}{\dim_H \mu}}}\text{ i.o.}\right\}<\dim_H \mu.\end{cases} $$

\end{corollary}

\section{Proof of the first item of Theorem \ref{Mainthm} }

We first start by recalling the following discrepency estimates for rotations on the  $2$-dimensional  Torus.
\begin{theoreme}[\cite{Saussogalato}, Proposition 17]
\label{thmDiscrep}
Let $\beta \in\mathbb{R}^2$ rational independent be such that $\gamma_\ell(\beta)<+\infty.$ Then for every $\gamma>\gamma_\ell (\beta),$ one has $$\sup\left\{\vert \frac{\sum_{0\leq k\leq n-1}\chi_R (k\beta)}{n}-\mathcal{L}^2(R)\vert, \ R\text{ is a rectangle }\right\}=O(\frac{1}{n^{\frac{1}{\gamma}}}) .$$
\end{theoreme}
Before proving the first item of Theorem \ref{Mainthm}, we start by a lemma, which easily deduced from the fact that  $\mathcal{L}^3$ is a doubling outer-regular measure. 
\begin{lemme}
Let  $\phi:\mathbb{N}\to \mathbb{R}_+$ be a non-increasing mapping. The following assertion are equivalent:
\begin{itemize}
\item[(1)] There exists $\gamma\geq 1,C>0$ such that the system $(S_{\alpha},\mathcal{L}^3)$ is $(\gamma,C\phi)$-mixing,\medskip
\item[(2)] there exists $\gamma\geq 1,C>0$ such that for every Borel set $A$ and every ball $B$, $$\mathcal{L}^3\Big(T^{-n}(A)\cap B\Big)\leq \gamma \mathcal{L}^3(A)\times \mathcal{L}^3(B)+C\phi(n)\mathcal{L}^3(A),$$
\item[(3)] there exists $\gamma \geq 1,C>0$ such that for every  dyadic cube $B$ and every dyadic cube $A$ with $\vert A \leq \vert B \vert$, $$\mathcal{L}^3\Big(T^{-n}(A)\cap B\Big)\leq \gamma \mathcal{L}^3(A)\times \mathcal{L}^3(B)+C\phi(n)\mathcal{L}^3(A).$$
\end{itemize}
\end{lemme}

We now establish the first item of Theorem \ref{Mainthm}. We fix $\gamma>\gamma_{\ell}=16,$ $A,B$ two dyadic cubes with $\vert A\vert \leq \vert B \vert.$ Let $n_B \in\mathbb{N}$ be the generation of $B$, i.e., the integer  such that $B\in\mathcal{D}_{n_B}.$ Let us write $$\begin{cases} A=I_A \times R_A  \\ B=I_B \times R_B ,\end{cases}$$  
where $I_A,I_B$ are  dyadic intervals and $R_A,R_B$ are dyadic squares.

In what follows, we will identify the a dyadic interval with its coding in base $2$. More precisely, let $f_0,f_1 :\mathbb{R}\to\mathbb{R}$ defined by $f_0(x)=\frac{x}{2}$ and $f_1(x)=\frac{x+1}{2}$ and $S=\left\{f_0,f_1\right\}.$ Denote by $\pi$ the canonical projection from $\left\{0,1\right\}^{\mathbb{N}}$ to $\mathbb{T}^1$, defined by $$\pi\Big((i_n)_{n\in\mathbb{N}}\Big)=\lim_{p\to +\infty}f_{i_1}\circ...\circ f_{i_p}(0)$$
and given a word $\underline{i}=(i_1,...,i_p)\in\left\{0,1\right\}^p,$ set $$[\underline{i}]=\left\{(i_1,...,i_p,(x_n)_{n\in\mathbb{N}}), \ (x_n)_{n\in\mathbb{N}}\in \left\{0,1\right\}^{\mathbb{N}}\right\}.$$

Given a dyadic interval $I\in\mathcal{D}_p$, we will write $I=[\underline{i}]$ where $\underline{i}\in\left\{0,1\right\}^p$ is such that $$I=\pi\Big([\underline{i}]\Big).$$
Let us write $I_A =[\underline{i}_A],I_B =[\underline{i}_B].$ Given $n\in\mathbb{N},$ one has

\begin{equation}
\label{EquaSmoinsn}
S_{\alpha}^{-n}(A)=\bigcup_{\underline{i}=(i_1,...,i_n)\in\left\{0,1\right\}^n}[\underline{i},\underline{i}_A]\times\Big(R_A-\alpha\#\left\{1\leq k\leq n: \ i_k=0\right\}\Big).
\end{equation}

So that $S_{\alpha}^{-n}(A)$ is a collection of $2^{n}$ disjoint $3$-dimensional rectangles which are products of $1$-dimensional intervals by $2$-dimensional squares of length $2^{-n}\mathcal{L}^1(I_A)\times \mathcal{L}^2(R_A).$ In addition, if $n\leq n_B,$ then at most one of these rectangles can intersect $B,$ so that, in this case,
\begin{align*}
\mathcal{L}^3\Big(S_{\alpha}^{-n}(A)\cap B\Big)&\leq 2^{-n}\mathcal{L}(I_A)\times \mathcal{L}^2(R_A)=2^{-n}\mathcal{L}^3(A)\leq \frac{1}{n^s}\mathcal{L}^3(A)\\
&\leq \frac{1}{n^s}\mathcal{L}^3(A)+\mathcal{L}^3(B)\times \mathcal{L}^3(A)
\end{align*}
for every $0< s\leq 1.$

We now assume that $n\geq n_B +1.$ Let $a,b \in\mathbb{T}^3$ be the centers of $A$ and $B$. Recall that $\vert R_A \vert \leq \vert R_B \vert,$ and notice that, if $1\leq k\leq n$ is such that $$\left\{ R_A-k\alpha \right\} \cap R_B \neq \emptyset,$$ one certainly has $$k\alpha\in \widetilde{R},$$ 
where $\widetilde{R}$ is a square centered on the center of $b-a$ with length-sides $3\vert R_B \vert.$ Thus, writing 
$$W_{n,A,B}=\left\{\underline{h}=(i_1,...,i_{n-n_B})\in\left\{0,1\right\}^{n-n_B}: \ \alpha\#\left\{1\leq k\leq n: \ i_k=0\right\}\in \widetilde{R} \right\} $$

 one has 
\begin{align*}
S_{\alpha}^{-n}(A)\cap B\subset\bigcup_{\underline{h}=(i_1,...,i_{n-n_B})\in W(n,A,B)}[\underline{i}_B,\underline{h},\underline{i}_A]\times\Big(R_A-\alpha\#\left\{1\leq k\leq n: \ i_k=0\right\}\Big).
\end{align*}  
This yields 
\begin{align}
\label{EquMesMix}
\mathcal{L}^3\Big(S_{\alpha}^{-n}(A)\cap B\Big)\leq \sum_{0\leq k\leq n-n_B}\chi_{\widetilde{R}}(k\alpha)\binom{n-n_B}{k}2^{-n}\mathcal{L}^1(I_A)\times \mathcal{L}^2(R_A).
\end{align}
We now estimate $\sum_{0\leq k\leq n-n_B}\chi_{\widetilde{R}}(k\alpha)\binom{n-n_B}{k}.$ Using Abel summation formula, one has 

\begin{align}
\label{EquaSigma1plusSigma2}
&\sum_{0\leq k\leq n-n_B}\chi_{\widetilde{R}}(k\alpha)\binom{n-n_B}{k}\nonumber\\
&=\sum_{0\leq k\leq n-n_B -1}\Big(\sum_{0\leq j\leq k}\chi_{\widetilde{R}}(j\alpha)\Big)\times \Big(\binom{n-n_B}{k}-\binom{n-n_B}{k+1}\Big)+\sum_{0\leq j\leq n-n_B}\chi_{\widetilde{R}}(j\alpha)\nonumber\\
&= \sum_{0\leq k\leq n-n_B -1}k\Big(\frac{\sum_{0\leq j\leq k}\chi_{\widetilde{R}}(j\alpha)}{k}-\mathcal{L}^2(\widetilde{R})\Big)\times \Big(\binom{n-n_B}{k}-\binom{n-n_B}{k+1}\Big)\nonumber\\
&+(n-n_B)\Big(\frac{\sum_{0\leq j\leq n-n_B}\chi_{\widetilde{R}}(j\alpha)}{n-n_B}-\mathcal{L}^2(\widetilde{R})\Big)+\nonumber\\
&\mathcal{L}^2(\widetilde{R})\Big(\sum_{0\leq k\leq n-n_B -1}k\times \Big(\binom{n-n_B}{k}-\binom{n-n_B}{k+1}\Big)+n-n_B\Big)=\Sigma_1 +\Sigma_2,
\end{align}
where 
\begin{align*}
\Sigma_1&= \sum_{0\leq k\leq n-n_B -1}k\Big(\frac{\sum_{0\leq j\leq k}\chi_{\widetilde{R}}(j\alpha)}{k}-\mathcal{L}^2(\widetilde{R})\Big)\times \Big(\binom{n-n_B}{k}-\binom{n-n_B}{k+1}\Big)\\
&+(n-n_B)\Big(\frac{\sum_{0\leq j\leq n-n_B}\chi_{\widetilde{R}}(j\alpha)}{n-n_B}-\mathcal{L}^2(\widetilde{R})\Big)\text{ and }\\
\Sigma_2&=\mathcal{L}^2(\widetilde{R})\Big(\sum_{0\leq k\leq n-n_B -1}k\times \Big(\binom{n-n_B}{k}-\binom{n-n_B}{k+1}\Big)+n-n_B\Big).
\end{align*}
Using again Abel summation formula on $\Sigma_2,$ one gets 
\begin{equation}
\label{estimsigma2}
\Sigma_2=\mathcal{L}^2(\widetilde{R})\sum_{0\leq k\leq n-n_B}\binom{n-n_B}{k} =\mathcal{L}^2(\widetilde{R})2^{n-n_B}.
\end{equation}

Recalling Theorem \ref{thmDiscrep}, there exists a constant $\kappa>0$ such that 
\begin{align*}
 &\sum_{0\leq k\leq n-n_B -1}k\Big(\frac{\sum_{0\leq j\leq k}\chi_{\widetilde{R}}(j\alpha)}{k}-\mathcal{L}^2(\widetilde{R})\Big)\times \Big(\binom{n-n_B}{k}-\binom{n-n_B}{k+1}\Big)\\
  &\leq \sum_{0\leq k\leq n-n_B}k\times \frac{\kappa}{k^{\frac{1}{\gamma}}}\Big \vert \binom{n-n_B}{k}-\binom{n-n_B}{k+1} \Big\vert\\
&\leq 2 \sum_{\frac{n-n_B}{2}\leq k\leq n-n_B}k\times \frac{\kappa}{k^{\frac{1}{\gamma}}}\Big( \binom{n-n_B}{k}-\binom{n-n_B}{k+1} \Big)\\
&\leq \frac{2^{1+\frac{1}{\gamma}}\kappa}{(n-n_B)^{\frac{1}{\gamma}}}\sum_{\frac{n-n_B}{2}\leq k\leq n-n_B}k\times\Big( \binom{n-n_B}{k}-\binom{n-n_B}{k+1} \Big)\\
&\leq \frac{3\times 2^{1+\frac{1}{\gamma}}\kappa}{(n-n_B)^{\frac{1}{\gamma}}}\sum_{0\leq k\leq n-n_B}k\times\Big( \binom{n-n_B}{k}-\binom{n-n_B}{k+1} \Big).
\end{align*}
Using the same argument on $(n-n_B)\Big(\frac{\sum_{0\leq j\leq n-n_B}\chi_{\widetilde{R}}(j\alpha)}{n-n_B}-\mathcal{L}^2(\widetilde{R})\Big),$ one obtains that 
\begin{align*}
\Sigma_1 \leq  \frac{3\times 2^{1+\frac{1}{\gamma}}\kappa}{(n-n_B)^{\frac{1}{\gamma}}}\Big(\sum_{0\leq k\leq n-n_B}k\times\Big( \binom{n-n_B}{k}-\binom{n-n_B}{k+1} \Big)+n-n_B\Big).
\end{align*}
and using Abel summation formula, one gets 
$$ \Sigma_1 \leq  \frac{3\times 2^{1+\frac{1}{\gamma}}\kappa}{(n-n_B)^{\frac{1}{\gamma}}}2^{n-n_B}.$$
Since $$\lim_{n_B \to+ \infty}\frac{2^{-n_B}}{(1-\frac{n_B}{n_B+1})^{\frac{1}{\gamma}}}=0,$$ there exists a constant $\eta>0$ such that 
\begin{equation}
\label{estisigma1}
\Sigma_1 \leq \eta \times\frac{ 2^{n}}{n^{\frac{1}{\gamma}}}.
\end{equation}

So, combining  \eqref{EquMesMix}, \eqref{EquaSigma1plusSigma2}, \eqref{estisigma1} and \eqref{estimsigma2}, one gets 
\begin{align*}
&\mathcal{L}^3\Big(S_{\alpha}^{-n}(A)\cap B\Big)\leq \\
& 2^{n-n_B}\mathcal{L}^2(\widetilde{R})\times  2^{-n}\times \mathcal{L}^1(I_A)\times \mathcal{L}^2(R_A)+ \eta \times\frac{ 2^{n}}{n^{\frac{1}{\gamma}}}\times  2^{-n}\times \mathcal{L}^1(I_A)\times \mathcal{L}^2(R_A)
\end{align*}
Recalling that $\mathcal{L}^{2}(\widetilde{R})=3^2 \mathcal{L}^2(R_B)$ and that $2^{-n_B}=\mathcal{L}^1(I_B),$ one gets $$\mathcal{L}^3\Big(S_{\alpha}^{-n}(A)\cap B\Big)\leq 9\mathcal{L}^3(A)\times \mathcal{L}^3(B)+\frac{\eta}{n^{\frac{1}{\gamma}}}\times\mathcal{L}^3(A).$$
Combining the case $n\leq n_B$ and $n\geq n_B+1,$ we conclude that $(S_{\alpha},\mathcal{L}^3)$ is $(9,n\mapsto\frac{\max\left\{\eta,1\right\}}{n^{\frac{1}{\gamma}}})$-mixing.

\section{Proof of the second item of Theorem \ref{Mainthm}}

We recall that, $\alpha=(\alpha_1,\alpha_2)$ being defined by Proposition \ref{alpha12}, one has

\begin{equation}
\label{equaqnqprimen}
\begin{cases}q_n^4 \leq q'_n \leq 4q_n^4 \\ {q'}_{n-1}^4\leq q_n \leq 4{q'}_{n-1}^4 \\
\gamma_{\ell}(\alpha_1,\alpha_2)=16,\end{cases}
\end{equation}

where $(\frac{p_n}{q_n})_{n\in\mathbb{N}}$ and  $(\frac{p'_n}{q'_n})_{n\in\mathbb{N}}$ are the partial quotient of $\alpha_1$ and $\alpha_2.$  Let us first notice that these equation implies in particular that

\begin{equation}
\label{Equaunplusksi}
\begin{cases} q_n^{16}\leq q_{n+1}\leq 16 q_n^{16} \\  {q'}_n^{16}\leq {q'}_{n+1}\leq 16 {q'}_n^{16}. \end{cases}
\end{equation}
We first show the following.
\begin{proposition}
\label{propoorbit0}
$$\dim_H \left\{y\in\mathbb{T}^2 : \ \vert\vert k\alpha-y\vert\vert_{\infty}\leq \frac{1}{k^{\frac{1}{3}}}\text{ i.o.}\right\}<2.$$

\end{proposition}

The following proof is inspired by the proof of \cite[Theorem 2]{bugchev}.
 \begin{proof}
Fix $\frac{12}{5}<\theta<\frac{46}{15}$ and $\frac{12}{5}<\beta<4$ and let us set, for every $n\in\mathbb{N},$ $Q_{2n}=q_n  q'_n$, $P_{2n}=Q_{2n}^{\theta},$ $Q_{2n+1}=q_{n+1}q'_n$ and $P_{2n+1}=Q_{2n+1}^{\beta}.$  By \eqref{equaqnqprimen}, one has $$\begin{cases}Q_{2n}\leq P_{2n}\leq Q_{2n+1} \\ Q_{2n+1}\leq P_{2n+1}\leq Q_{2n+2}.\end{cases} $$

Notice in addition that $$\begin{cases} \bigcup_{Q_{2n}\leq k\leq Q_{2n+1} }B(k\alpha,\frac{1}{k^{\frac{1}{3}}}) \subset \bigcup_{Q_{2n}\leq k\leq P_{2n} }B(k\alpha,\frac{1}{Q_{2n}^{\frac{1}{3}}})\bigcup \bigcup_{P_{2n}+1 \leq k\leq Q_{2n+1}}B(k\alpha,\frac{1}{P_{2n}^{\frac{1}{3}}})\\
\bigcup_{Q_{2n+1}\leq k\leq Q_{2n+2} }B(k\alpha,\frac{1}{k^{\frac{1}{3}}}) \subset \bigcup_{Q_{2n+1}\leq k\leq P_{2n+1} }B(k\alpha,\frac{1}{Q_{2n+1}^{\frac{1}{3}}})\bigcup \bigcup_{P_{2n+1}+1 \leq k\leq Q_{2n+2}}B(k\alpha,\frac{1}{P_{2n+1}^{\frac{1}{3}}}).\end{cases} $$

In addition, since $\left\{Q_{2n}\leq k< Q_{2n+2}\right\}$ is a partition of $\mathbb{N},$ one has 
\begin{align*}
&\limsup_{k\to+\infty}B(k\alpha,\frac{1}{k^{\frac{1}{3}}})=\limsup_{n\to+\infty}\bigcup_{Q_{2n}\leq k\leq Q_{2n+1} }B(k\alpha,\frac{1}{k^{\frac{1}{3}}})\bigcup  \bigcup_{Q_{2n+1}\leq k\leq Q_{2n+2} }B(k\alpha,\frac{1}{k^{\frac{1}{3}}}) \\
&\subset \limsup_{n\to+\infty}\bigcup_{Q_{2n}\leq k\leq P_{2n} }B(k\alpha,\frac{1}{Q_{2n}^{\frac{1}{3}}}) \bigcup \limsup_{n\to+\infty}\bigcup_{P_{2n}+1 \leq k\leq Q_{2n+1}}B(k\alpha,\frac{1}{P_{2n}^{\frac{1}{3}}})\bigcup\\
&\limsup_{n\to+\infty}\bigcup_{Q_{2n+1}\leq k\leq P_{2n+1} }B(k\alpha,\frac{1}{Q_{2n+1}^{\frac{1}{3}}})\bigcup \limsup_{n\to+\infty}\bigcup_{P_{2n+1}+1 \leq k\leq Q_{2n+2}}B(k\alpha,\frac{1}{P_{2n+1}^{\frac{1}{3}}}).
\end{align*}

On what follows, we provide an upper-bound for $$\dim_H \limsup_{n\to+\infty}\bigcup_{Q_{2n}\leq k\leq P_{2n} }B(k\alpha,\frac{1}{Q_{2n}^{\frac{1}{3}}}) \text{ and }\dim_H \limsup_{n\to+\infty}\bigcup_{P_{2n}+1 \leq k\leq Q_{2n+1}}B(k\alpha,\frac{1}{P_{2n}^{\frac{1}{3}}}).$$ Upper-bounds for $$\dim_H \limsup_{n\to+\infty}\bigcup_{Q_{2n+1}\leq k\leq P_{2n+1} }B(k\alpha,\frac{1}{Q_{2n+1}^{\frac{1}{3}}})\text{ and }\dim_H \limsup_{n\to+\infty}\bigcup_{P_{2n+1}+1 \leq k\leq Q_{2n+2}}B(k\alpha,\frac{1}{P_{2n+1}^{\frac{1}{3}}})$$ are obtained using similar arguments.

Observe that $$\max\left\{\vert \alpha_1 -\frac{p_n}{q_n}\vert,\vert \alpha_2 -\frac{p'_n}{q'_n}\vert\right\}\leq \max\left\{\frac{1}{q_n q_{n+1}},\frac{1}{q'_n q'_{n+1}}\right\}\leq \frac{1}{q_n q_{n+1}}.$$ Thus, one has, for large enough $n$, 
\begin{align*}
P_{2n}\vert\vert\alpha -(\frac{p_n}{q_n},\frac{p'_n}{q'_n}) \vert\vert_{\infty}&\leq \Big(q'_n q_n\Big)^{\theta}\times \frac{1}{q_n q_{n+1}}\leq 4^{\frac{5}{2}}q_n^{5\theta-17}
\end{align*}
and $Q_{2n}^{-\frac{1}{3}}\leq q_n^{-\frac{5}{3}}.$
Since $17 -5\theta> \frac{5}{3},$ one has $$P_{2n}\vert\vert\alpha -(\frac{p_n}{q_n},\frac{p'_n}{q'_n}) \vert\vert_{\infty}\leq Q_{2n}^{-\frac{1}{3}}.$$
This implies that 
\begin{align*}
\bigcup_{Q_{2n}\leq k\leq P_{2n} }B(k\alpha,\frac{1}{Q_{2n}^{\frac{1}{3}}})\subset \bigcup_{0\leq k\leq P_{2n}}B\Big(k(\frac{p_n}{q_n},\frac{p'_n}{q'_n}),\frac{2}{Q_{2n}^{\frac{1}{3}}}\Big).
\end{align*}
Given $0\leq k\leq q_n,$ write $\ell_{k}=\left\{(\frac{k}{q_n},t), \ t\in[0,1]\right\}$ and $$U(\ell_k,\frac{2}{Q_{2n}^{\frac{1}{3}}})=\bigcup_{x\in \ell_k}B(x,\frac{2}{Q_{2n}^{\frac{1}{3}}}).$$
Clearly, $$\bigcup_{0\leq k\leq P_{2n}}B\Big(k(\frac{p_n}{q_n},\frac{p'_n}{q'_n}),\frac{2}{Q_{2n}^{\frac{1}{3}}}\Big)\subset \bigcup_{0\leq k\leq q_n}U(\ell_k,\frac{2}{Q_{2n}^{\frac{1}{3}}}).$$
Now, fix $0\leq s\leq 2,$ and let $\mathcal{C}_{\ell_k}$ be a  covering of each line $\ell_k$ by $Q_{2n}^{\frac{1}{3}}$ balls of radius $\frac{2}{Q_{2n}^{\frac{1}{3}}},$ one has 
\begin{align*}
\sum_{B\in \bigcup_{0\leq k\leq q_n}\mathcal{C}_{\ell_k}}\vert B\vert^s \leq \frac{2q_n}{Q_{2n}^{\frac{1}{3}}}\times \frac{2}{Q_{2n}^{\frac{s}{3}}}\leq Cq_n Q_{2n}^{-\frac{s-1}{3}}\leq Cq_n^{1-\frac{5\theta}{3}(s-1)}.
\end{align*}
since $\theta>\frac{3}{5},$ one has $1-\frac{5\theta}{3}<0$, so that there exists $1<s_1<2$ such that $1-\frac{5\theta}{3}(s_1-1)<0.$ Since $q_n\geq q_1 ^{16n},$ one has $$\sum_{n\geq 1}Cq_{n}^{1-\frac{5\theta}{3}(s_1-1)}<+\infty.$$
This proves that $$\dim_H \limsup_{n\to+\infty}\bigcup_{Q_{2n}\leq k\leq P_{2n} }B(k\alpha,\frac{1}{Q_{2n}^{\frac{1}{3}}})\leq s_1<2.$$

We now show that $$\dim_H \limsup_{n\to+\infty}\bigcup_{P_{2n}+1 \leq k\leq Q_{2n+1}}B(k\alpha,\frac{1}{P_{2n}^{\frac{1}{3}}})<2.$$
Similarly as before,
\begin{align*}
Q_{2n+1}\vert\vert\alpha-(\frac{p_{n+1}}{q_{n+1}},\frac{p'_n}{q'_n}) \vert\vert_{\infty}\leq q_n^{16+4}\times \frac{C}{q'_n q'_{n+1}}\leq C q_n^{20-(4+4\times 16)}=C q_n^{-48}.
\end{align*} 
 And $P_n^{-\frac{1}{3}}\leq q_n^{\frac{-5\theta}{3}} .$ Since $48>\frac{5\theta}{3},$ one has $$Q_{2n+1}\vert\vert\alpha-(\frac{p_{n+1}}{q_{n+1}},\frac{p'_n}{q'_n}) \vert\vert_{\infty}\leq P_n^{-\frac{1}{3}},$$
 which implies that
 \begin{align*}
\bigcup_{P_{2n}\leq k\leq Q_{2n+1} }B(k\alpha,\frac{1}{P_{2n}^{\frac{1}{3}}})\subset \bigcup_{0\leq k\leq Q_{2n+1}}B\Big(k(\frac{p_{n+1}}{q_{n+1}},\frac{p'_n}{q'_n}),\frac{2}{P_{2n}^{\frac{1}{3}}}\Big).
\end{align*}
Thus, using the same argument as before, there exists $q'_n +1$ lines of length $1$ containing $\left\{k(\frac{p_{n+1}}{q_{n+1}},\frac{p'_n}{q'_n})\right\}_{0\leq k\leq Q_{2n+1}}.$ Thus by covering each of these lines by balls of length-sides $\frac{2}{P_{2n}^{\frac{1}{3}}},$ one obtain a covering $\mathcal{C}$ of  $\bigcup_{P_{2n}\leq k\leq Q_{2n+1} }B(k\alpha,\frac{1}{P_{2n}^{\frac{1}{3}}})$ such that, for every $1<s<2,$ 
\begin{align*}
\sum_{B\in \mathcal{C}}\vert B\vert^s \leq q'_n P_{2n}^{-\frac{(s-1)}{3}}\leq Cq_n^{4-5\theta\frac{s-1}{3}}.
\end{align*}
Since $4-5\frac{\theta}{3}<0,$ there exists $1<s_2<2$ such that $4-5\theta\frac{s_2-1}{3}<0.$ Since $$\sum_{n\geq 1}Cq_n^{4-5\theta\frac{s_2-1}{3}}<+\infty, $$
one has $$\dim_H \limsup_{n\to+\infty}\bigcup_{P_{2n}+1 \leq k\leq Q_{2n+1}}B(k\alpha,\frac{1}{P_{2n}^{\frac{1}{3}}})\leq s_2<2,$$
which concludes the proof.
\end{proof}

Fixing any $t\in\mathbb{T}^2,$ since $$\vert\vert t+k\alpha -y\vert\vert_{\infty}\leq \frac{1}{k^{\frac{1}{3}}}\Leftrightarrow \vert\vert k\alpha -(y-t)\vert\vert_{\infty}\leq \frac{1}{k^{\frac{1}{3}}},$$
Proposition \ref{propoorbit0} implies that for every $t\in\mathbb{T}^2,$ one has $$\dim_H \left\{y\in\mathbb{T}^2 : \ \vert\vert t+k\alpha-y\vert\vert_{\infty}\leq \frac{1}{k^{\frac{1}{3}}}\text{ i.o.}\right\}<2.$$
Moreover, for ever $x\in\mathbb{T}^1,t\in\mathbb{T}^2,$
\begin{align*}
&\left\{y=(y_1,y_2,y_3)\in\mathbb{T}^3 : \ \vert\vert S_{\alpha}^n(x,t)-y \vert\vert_{\infty}\leq \frac{1}{n^{\frac{1}{3}}}\text{ i.o. }\right\}\\
&\subset \mathbb{T}^1 \times \left\{y'=(y_2,y_3)\in\mathbb{T}^2 : \  \vert\vert t+n\alpha-y'\vert\vert_{\infty}\leq \frac{1}{n^{\frac{1}{3}}}\text{ i.o.}\right\},
\end{align*}
which yields

\begin{align*}
&\dim_H \left\{y\in\mathbb{T}^3 : \ \vert\vert S_{\alpha}^n(x,t)-y \vert\vert_{\infty}\leq \frac{1}{n^{\frac{1}{3}}}\text{ i.o. }\right\} \\
&\leq 1+\dim_H \left\{y'\in\mathbb{T}^2 : \  \vert\vert t+n\alpha-y'\vert\vert_{\infty}\leq \frac{1}{n^{\frac{1}{3}}}\text{ i.o.}\right\}<3.
\end{align*}
 This proves Item $(2)$ of Theorem \ref{Mainthm}. 
  
\bibliographystyle{plain}
\bibliography{bibliogenubi}

\end{document}